\documentclass{article}
\pdfoutput=1

\usepackage[narrow]{arxiv}

\usepackage[utf8]{inputenc} % allow utf-8 input
\usepackage[T1]{fontenc}    % use 8-bit T1 fonts
\usepackage[hypertexnames=false]{hyperref}       % hyperlinks
\usepackage{url}            % simple URL typesetting
\usepackage{booktabs}       % professional-quality tables
\usepackage{amsfonts}       % blackboard math symbols
\usepackage{nicefrac}       % compact symbols for 1/2, etc.
\usepackage{microtype}      % microtypography

\usepackage{lipsum}         % Can be removed after putting your text content
\usepackage{graphicx}
\usepackage{natbib}
\usepackage{doi}

\title{Near-Optimal Methods for Generalized Lipschitz Convex Optimization}

\newif\ifuniqueAffiliation
\uniqueAffiliationtrue

\ifuniqueAffiliation % Standard variant of author block
\author{
  Dmitry Kovalev \\
  Yandex Research \\
  \texttt{dakovalev1@gmail.com}
  \And
  Grigory Malinovsky\\
  KAUST\thanks{King Abdullah University of Science and Technology}\\
  \texttt{g.s.malinovsky@gmail.com}
}
\else
\usepackage{authblk}

\author[1]{%
  Dmitry Kovalev\thanks{\texttt{dakovalev1@gmail.com}}%
}
\affil[1]{Yandex Research}
\fi
\usepackage{algorithm,algpseudocodex}
\usepackage{nicematrix}

\usepackage[textsize=tiny,color=orange!20,bordercolor=orange,]{todonotes}
\presetkeys{todonotes}{inline}{}

\usepackage{xcolor}
\hypersetup{
  colorlinks=true,
  linkcolor=red!75!black,
  citecolor=blue!50!black
}

\usepackage{mathtools,amsthm,amssymb}
\allowdisplaybreaks[4]

\usepackage{enumitem}
\setlist{topsep=0pt,partopsep=0pt,itemsep=0pt}

\makeatletter\AddToHook{env/algorithmic/before}{\def\@currentcounter{ALG@line}}\makeatother
\usepackage{cleveref}
\crefalias{ALG@line}{line}
\usepackage{annotate}
\usepackage{mythm}
\usepackage{mymath}

\newcommand{\lip}{\ell^2}
\newcommand{\lipsqrt}{\ell}

\newcommand{\Rd}{\R^d}
\newcommand{\N}{\mathbb{N}}
\newcommand{\be}{\mathbf{e}}
\newcommand{\clip}{\operatorname{clip}\brr}
\newcommand{\hDelta}{\hat{\Delta}}

\begin{document}
\maketitle

\begin{abstract}
 In this paper, we consider a convex optimization problem, where the upper bound on the subgradient norm depends on the objective function. We establish a lower bound on the number of iterations required by any first-order method to find an approximate solution to the problem. Furthermore, we present two near-optimal algorithms that can match the lower bound up to logarithmic factors. The first algorithm is a projected GD with a carefully designed stepsize schedule which incorporates warm-up and decay phases. The second algorithm is scalar-stepsize clipped AdamW which extends the near-optimal guarantees to the stochastic setting, where only a noisy unbiased estimator of the subgradient is available.
\end{abstract}

\section{Introduction}

In this paper, we consider the problem of minimizing a convex and possibly non-smooth objective function $f\brr{}\colon \Rd \to \R$,
\begin{equation}\label{eq:main}
  \min_{x \in \Rd} f(x).
\end{equation}
We assume that at least one global solution $x^* \in \Rd$ to the problem exists and denote the corresponding optimal objective function value as $f^* = f(x^*)$.

We also define the mapping $\nabla f\brr{} \colon \Rd \to \Rd$ as an arbitrary selection from the subdifferential mapping of the objective function $f\brr{}$. By standard results in convex analysis \citep{rockafellar1997convex}, this mapping is well-defined and possibly non-unique, and, with a slight abuse of language, we refer to $\nabla f(x)$ as the gradient of the objective function $f\brr{}$ at the point $x \in \Rd$.

\subsection{Generalized Lipschitz Functions}

In this paper, we focus on the following key assumption on the gradient of the objective function:
\begin{aequation}<ass:Lip>
  \sqn{\nabla f(x)} \leq \lip([f(x) - f^*]^2)
  \quad\text{for all}\;\;
  x \in \Rd,
\end{aequation}
where $\lip\brr{}\colon \R_+ \to \R_{++}$ is a positive non-decreasing continuous concave function. When $\lip\brr{}$ is a constant function, this assumption is equivalent to the Lipschitz continuity of the objective function $f\brr{}$. Consequently, we will refer to \Cref{ass:Lip} as the generalized Lipschitzness. This assumption was previously studied by \citet[Definition~4]{grimmer2017convergence}. It also extends several popular existing assumptions, including $(\cL_0,\cL_1)$-smoothness \citep{zhang2019gradient}, $(\cH_0,\cH_1)$-smoothness \citep{alimisis2025we}, and $(\cM_0,\cM_1)$-Lipschitzness \citep{kovalev2026stochastic}. We also give an equivalent characterization of \Cref{ass:Lip} in the following \Cref{lem:lip}.
\begin{lemma}<lem:lip>
  \Cref{ass:Lip} holds if and only if the function $x \mapsto \Phi(f(x) - f^*)$ is $1$-Lipschitz, where the function $\Phi\brr{}\colon \R_+ \to \R_+$ is defined as follows:
  \begin{equation}\label{eq:Phi}
     \Phi(s) = {\textstyle\int_0^s} \mathrm{d}t/\lipsqrt(t^2),
     \quad\text{where}\;\;
     \lipsqrt\brr{t} = \sqrt{\lip\brr{t}}.
  \end{equation}
\end{lemma}

In this paper, we are also interested in the stochastic optimization setting, where we have access to the unbiased noisy estimator $\nabla_\xi f(x)$ of the gradient $\nabla f(x)$ at every point $x \in \Rd$, where $\xi \sim \cD$ is a random variable sampled from the distribution $\cD$. Consequently, we consider the following stochastic extension of \Cref{ass:Lip}:
\begin{aequation}<ass:Lip_stoch>
  \E[\xi\sim\cD]{\nabla_\xi f(x)} = \nabla f(x)
  ,\quad
  \E[\xi\sim\cD]{\sqn{\nabla_\xi f(x)}} \leq \lip([f(x) - f^*]^2)
  \quad\text{for all}\;\;
  x \in \Rd,
\end{aequation}
where $\lip\brr{}\colon \R_+ \to \R_{++}$ is a positive non-decreasing continuous concave function.
This assumption obviously implies its deterministic variant, \Cref{ass:Lip}, and it was also previously studied by \citet[Definition~4]{grimmer2017convergence}.

\subsection{Summary of Contributions}

In this paper, we establish the worst-case iteration complexity of first-order methods on generalized Lipschitz convex optimization problems in both deterministic and stochastic settings. More concretely, we make the following contributions:

\begin{enumerate}[label=\bf(\roman*)]
  \item In \Cref{sec:lb}, we establish the lower bound on the number of iterations required by any first-order optimization algorithm to find an $\epsilon$-approximate solution\footnote{A point $x \in \Rd$ is called an $\epsilon$-accurate solution to problem~\eqref{eq:main}, if $f(x) - f^* \leq \epsilon$.} to problem~\eqref{eq:main} satisfying \Cref{ass:Lip,ass:Lip_stoch}.
  \item In \Cref{sec:det}, we show that the standard projected gradient descent with scheduled stepsizes is a near-optimal algorithm for solving problem~\eqref{eq:main} under \Cref{ass:Lip}, i.e., it matches the lower complexity bound in \Cref{sec:lb} up to logarithmic factors.
  \item In \Cref{sec:stoch}, we show that scalar-stepsize clipped AdamW is a near-optimal algorithm for solving the stochastic version of problem~\eqref{eq:main} under \Cref{ass:Lip_stoch}, i.e., it also matches the lower bound in \Cref{sec:lb} up to logarithmic factors.
\end{enumerate}

\textbf{Notation.} In this paper, we use the following notations:
$\norm{}$ denotes the standard Euclidean norm induced by the standard scalar product $\<,>$ in $\Rd$;
for arbitrary vector $x \in \Rd$, $x_{[i]}$ denotes the $i$-th coordinate of $x$;
$\clip{}\colon v \mapsto \min\brf{1,1/\norm{v}}v$ is the standard clipping operation. We also use the convention $0/0 = 0$.

\section{Lower Complexity Bounds}\label{sec:lb}

In this section, we develop a lower bound on the number of iterations required by any first-order optimization algorithm to find an $\epsilon$-accurate solution to convex optimization problem~\eqref{eq:main} with a generalized Lipschitz convex objective function $f\brr{}$. For simplicity, we focus on the lower complexity bound for the deterministic problem satisfying \Cref{ass:Lip}, and this bound will obviously apply to the stochastic problem satisfying \Cref{ass:Lip_stoch}. We first specify the class of algorithms for which the lower bound applies. More specifically, we assume that the algorithm starts with the zero point $x_0 = 0 \in \Rd$ and generates points $x_k \in \Rd$ for $k \in \N$ satisfying the following condition:
\begin{aequation}<ass:span>
  x_k \in \spanset{\brf{\nabla f(x_0),\dots,\nabla f(x_{k-1})}}.
\end{aequation}
\Cref{ass:span} is referred to in literature as the {\em linear span assumption}. It is commonly used for constructing lower complexity bounds for various optimization problem classes \citep{nesterov2013introductory} because most practical algorithms satisfy this assumption. This assumption could be removed using slightly more complex reasoning (see e.g.~\citet{arjevani2019oracle}); however, such extension would not provide any new theoretical insights and would only introduce unnecessarily technical complications.

Our lower complexity bound is based on a specific hard instance of problem~\eqref{eq:main}. First, we define function $h_{d,\rho}\brr{}\colon \Rd \to \R_+$ as follows:
\begin{equation}\label{eq:h_hard}
  h_{d,\rho}(x) = \rho + \max\brf{x_{[1]},\dots,x_{[d]},-\rho}.
\end{equation}
This function is used to prove lower complexity bounds for standard Lipschitz convex optimization \citep{nesterov2013introductory}. 
To establish the lower bounds under the generalized Lipschitz assumption, we consider the objective function $f\brr{}$ of the following form:
\begin{equation}\label{eq:f_hard}
  f(x) = \cL(h_{d,\rho}(x))
  ,\quad\text{where}\;\;
  \cL(s) = \Phi^{-1}(s).
\end{equation}
In the following \Cref{lem:f_hard} we show that the function $f\brr{}$ in \cref{eq:f_hard} is well-defined and establish its basic properties.
\begin{lemma}<lem:f_hard>
  The function $f\brr{}$ in \cref{eq:f_hard} is well-defined, convex and satisfies \Cref{ass:Lip}. Moreover, it holds that $\partial f(x) \supset \lipsqrt([f(x) - f^*]^2)\cdot\partial h_{d,\rho}(x)$ for all $x \in \Rd$, and there exists a global solution $x^* \in \Rd$ to problem~\eqref{eq:main} satisfying $\norm{x^*} \leq \sqrt{d}\rho$.
\end{lemma}

Next, we choose the gradient $\nabla f\brr{}$ for the function $f\brr{}$ defined in \cref{eq:f_hard} as follows:
\begin{equation}\label{eq:grad_hard}
  \nabla f(x) = \lipsqrt([f(x) - f^*]^2)\nabla h_{d,\rho}(x)
  ,\;\;\text{where}\;\;
  \nabla h_{d,\rho}(x) = \begin{cases}
    \be_{i(x)} &  h_{d,\rho}(x) > 0\\
    0 & h_{d,\rho}(x) = 0
  \end{cases},
\end{equation}
and where $\be_i \in \Rd$ is the $i$-th standard unit basis vector, and $i(x)$ is defined as follows:
\begin{equation}
  i(x) = \min\brf{j \in \brf{1,\dots,d} : \rho + x_{[j]} = h_{d,\rho}(x)}.
\end{equation}
One can verify that $\nabla h_{d,\rho}(x) \in \partial h_{d,\rho}(x)$, and \Cref{lem:f_hard} implies that $\nabla f(x) \in \partial f(x)$ for all $x \in \Rd$.
In the following \Cref{lem:grad_hard}, we establish an important conclusion about the optimization trajectory of any optimization algorithm satisfying \Cref{ass:span} applied to our hard problem instance with the gradient $\nabla f \brr{}$ defined in \cref{eq:grad_hard}. This lemma is a trivial extension of the arguments by \citet{nesterov2013introductory} so its proof is omitted.
\begin{lemma}<lem:grad_hard>
  Let the gradient $\nabla f\brr{}$ be defined in \cref{eq:grad_hard}, and suppose that \Cref{ass:span} holds. Then $x_k \in \spanset{\brf{\be_1,\dots,\be_k}}$ for all $k \in \brf{1,\dots,d}$
\end{lemma}
Now, we are ready to establish the desired lower complexity bound in the following \Cref{thm:lb}.
\begin{theorem}<thm:lb>
  Let $\cR > 0$ and $\epsilon\in(0,\Phi^{-1}(\cR)]$ be arbitrary positive constants. Let $\lip\brr{}\colon \R_+ \to \R_{++}$ be an arbitrary positive non-decreasing continuous concave function. Then, there exist a convex function $f\brr{}\colon \Rd \to \R$ and a corresponding gradient $\nabla f\brr{}\colon \Rd \to \Rd$ satisfying \Cref{ass:Lip}, such that there exists a solution $x^* \in \Rd$ to problem~\eqref{eq:main} satisfying $\norm{x^*} \leq \cR$, and, to output a point $x_K \in \Rd$ satisfying $f(x_K) - f^* < \epsilon$, any algorithm satisfying \Cref{ass:span} requires at least the following number of iterations:
\begin{equation}
  K = \tilde\Omega\brr*{[\cR\lipsqrt(\epsilon^2)/\epsilon]^2},
\end{equation}
where $\tilde\Omega\brr{}$ hides logarithmic factors depending on $\lipsqrt(0)$ and $\lipsqrt(\epsilon^2)$.
\end{theorem}

We make two remarks regarding the result in \Cref{thm:lb}. First, the assumption $\epsilon \leq \Phi^{-1}(\cR)$ is not restrictive, because under \Cref{ass:Lip} it holds that $f(x_0) - f^* \leq \Phi^{-1}(\cR)$ due to \Cref{lem:lip}. Second, in the case where the function $\lip\brr{}$ is constant, \Cref{ass:Lip} is equivalent to the $\lipsqrt(0)$-Lipschitzness of the function $f(x)$, \citet{nesterov2013introductory} gives the standard lower complexity bound $\Omega\brr*{[\cR\lipsqrt(0)/\epsilon]^2}$. Consequently, \Cref{thm:lb} extends this result to \Cref{ass:Lip} with arbitrary positive non-decreasing continuous concave function $\lip\brr{}$.

\section{Near-Optimal Algorithm for Deterministic Problems}\label{sec:det}

In this section, we consider a simple projected gradient descent for solving problem~\eqref{eq:main} with the following iterations:
\begin{equation}\label{eq:GD}
  x_0 = 0,\qquad x_{k+1} = \proj[Q_\cR]{x_k - \eta_k \nabla f(x_k)},
\end{equation}
where $\eta_k > 0$ are positive stepsizes and $Q_\cR = \brf{x \in \Rd : \norm{x} \leq \cR}$ is a ball of radius $\cR \geq \norm{x^*}$ centered at the origin. Using similar arguments to \citet[Appendix~A]{kovalev2026stochastic}, one can show that this algorithm may not reach lower complexity bound in \Cref{thm:lb} even up to logarithmic factors, when the stepsizes $\eta_k$ are constant. However, we will show that with a carefully designed stepsize schedule, projected gradient descent can reach this bound up to logarithmic factors.

We start with the following \Cref{lem:GD}, which establishes an essential inequality used in the last-iterate convergence analysis of gradient descent. This lemma is a consequence of Lemma~2.1 of \citet{zamani2025exact}. We present its simplified proof for completeness.
\begin{lemma}<lem:GD>
  The iterations in \cref{eq:GD} satisfy the following inequality for all $x \in Q_\cR$ and $0 \leq l \leq k$:
  \begin{equation}
    f(x_k) - f(x) \leq \tmfrac{\sqn{x_l - x}}{4H_{k,l}} + \msum_{j=l}^{k-1} \tmfrac{\eta_j^2\sqn{\nabla f(x_j)}}{4H_{k,j+1}} + \eta_k\sqn{\nabla f(x_k)},
    \quad\text{where}\;\;
    H_{k,j} = \msum_{i=j}^k \eta_i.
  \end{equation}
\end{lemma}
Next, we specify the stepsize schedule used in \cref{eq:GD}. Given a positive constant $F_0 > 0$, we define the stepsizes $\eta_k$ recursively as follows:
\begin{equation}\label{eq:stepsize}
  \eta_k = \tmfrac{F_k}{L_k\lip(F_k^2)},\quad F_{k+1} = (1-\alpha_k)F_k,
  \quad\text{where}\;\;
  L_k = 64 + 16\ln(k+2)\;\;\text{and}\;\;
  \alpha_k = \tmfrac{\eta_kF_k}{8\cR^2}.
\end{equation}
Using this stepsize schedule, we can establish the key \Cref{thm:GD}, which establish an upper bound on the objective function at the points generated by \cref{eq:GD,eq:stepsize}. We also establish an explicit iteration complexity in \Cref{cor:GD}.
\begin{theorem}<thm:GD>
  Let \Cref{ass:Lip} hold. Let $\Delta_0 \geq f(x_0) - f^* > 0$ and $\Delta_0 \leq \cR \lipsqrt(\Delta_0^2)$, and let $F_0 = 3\Delta_0$. Then the iterates in \cref{eq:GD} with the stepsize schedule in \cref{eq:stepsize} satisfy $f(x_k) - f^* \leq F_k$.
\end{theorem}
\begin{corollary}<cor:GD>
  Under the conditions of \Cref{thm:GD}, suppose that $\epsilon > 0$. Then, it is sufficient to perform the following number of iterations $K$ to ensure $f(x_K) - f^* \leq \epsilon$:
  \begin{equation}
    K = \tilde{\cO}\brr*{\brs{\cR\lipsqrt(\epsilon^2)/\epsilon}^2},
  \end{equation}
  where $\tilde\cO\brr{}$ hides logarithmic factors depending on $\lipsqrt(\epsilon^2)$, $\cR$, $F_0$, and $\epsilon$.
\end{corollary}
One can show that the stepsize schedule in \cref{eq:stepsize} may incorporate an increasing stage followed by a decay stage, for instance, in the case where $\lip\brr{}$ is an affine function. \Cref{thm:lb} and \Cref{cor:GD} show that it achieves near optimal convergence guarantees. Consequently, we provide theoretical support for the common practice of specifying stepsize schedules with warm-up and decay in advance. 
There are several existing works \citep{liu2025theoretical,alimisis2025we,riabinin2026does} that claim to analyze gradient methods with stepsize warm-up under various special cases of \Cref{ass:Lip}. However, these works do not analyze predefined stepsize schedules, but rather analyze specific adaptive stepsizes depending on the current objective gap. Such stepsizes are hard to implement, because they require $f^*$, and, to the best of our knowledge, have seen little adoption in practice. Consequently, in contrast to our results, these works leave open the question of whether predefined stepsize schedules can offer provable benefits.

\section{Near-Optimal Algorithm for Stochastic Problems}\label{sec:stoch}

\begin{algorithm}[t]
  \caption{Scalar-stepsize Clipped AdamW}
  \label{alg:adamw}
  \begin{algorithmic}[1]
    \State \textbf{input:} $\alpha\in (0,1)$, $\cR > \norm{x^*}$, $K \in \N$
    \State $x_0 = m_0 = 0 \in \Rd$, $s_0 = 0$
    \For{$k = 0,\dots,K-1$}
    \State $x_{k+1} = (1-\alpha)x_k - \alpha\cR\clip{m_k / \sqrt{s_k}}$
    \label{line:x}
    \State $g_{k+1} = \nabla_{\xi_{k+1}} f(x_{k+1})$, where $\xi_{k+1} \sim \cD$ is a fresh sample
    \label{line:g}
    \State $m_{k+1} = (1-\alpha)m_k + \alpha g_{k+1}$,\;\; $s_{k+1} = (1-\alpha)^2 s_k + \alpha^2 \sqn{g_{k+1}}$\label{line:ms}
    \EndFor
    \State \textbf{output:} $x_K \in \Rd$
  \end{algorithmic}
\end{algorithm}

 Unfortunately, we are able to establish the near-optimal convergence guarantees for projected gradient descent with scheduled stepsize in \Cref{sec:det} only for deterministic problems under \Cref{ass:Lip}. However, in this section, we will show that a scalar-stepsize version of clipped AdamW can achieve near-optimal guarantees for stochastic problems under \Cref{ass:Lip_stoch}. The update rule for this algorithm is given in \Cref{alg:adamw}.

 Further, we establish the convergence analysis for \Cref{alg:adamw}. It is inspired by the existing analysis of \citet{kovalev2026stochastic} and is substantially simplified and improved. We start with \Cref{lem:adamw}, which establishes the key inequality used in the analysis of \Cref{alg:adamw}.
 \begin{lemma}<lem:adamw>
  Let \Cref{ass:Lip_stoch} hold, and let $\Delta_k^2 = \E{(f(x_k) - f^*)^2}$. Then, the iterations of \Cref{alg:adamw} satisfy the following inequality:
  \begin{equation}
    \Delta_K^2 \leq 3(1-\alpha)^{2K}\Delta_0^2 + [8\alpha\cR]^2\msum_{k=1}^K (1-\alpha)^{2(K-k)}\lip(\Delta_k^2)
  \end{equation}
 \end{lemma}
 Next, we show that the function $\lipsqrt\brr{}$ can be upper-bounded with a certain affine function. In particular, we show that \Cref{ass:Lip,ass:Lip_stoch} imply the $(\cM_0,\cM_1)$-Lipschitzness of the function $f\brr{}$. Moreover, we will further show, that we do not lose much information with this reduction and we will still be able to obtain near optimal convergence guarantees under the general \Cref{ass:Lip_stoch}.
 \begin{lemma}<lem:M01>
  Let $\cM_0 = \lipsqrt(\epsilon^2)$ and $\cM_1 = \lipsqrt(\epsilon^2)/\epsilon$ for an arbitrary $\epsilon > 0$. Then, for all $s \geq 0$,
  \begin{equation}
    \lip(s^2) \leq \cM_0^2 + \cM_1^2 s^2.
  \end{equation}
 \end{lemma}

  Finally, we establish the convergence guarantees for \Cref{alg:adamw}. The following \Cref{thm:adamw} shows that the iteration complexity of \Cref{alg:adamw} matches the lower complexity bounds established in \Cref{thm:lb}. Consequently, our result shows that scalar-stepsize clipped AdamW is a near optimal algorithm for solving stochastic convex optimization problems under \Cref{ass:Lip_stoch}.

 \begin{theorem}<thm:adamw> Suppose that \Cref{ass:Lip_stoch} holds. Let $\alpha = 1/\max\brf{2, [16\cR\lipsqrt(\epsilon^2)/\epsilon]^2}$, where $\epsilon \in (0, f(x_0) - f^*)$. Then, to reach the precision $\E{f(x_K) - f^*} \leq \epsilon$, it is sufficient to perform the following number of iterations of \Cref{alg:adamw}:
  \begin{equation}
    K = \ceil*{\tmfrac{1}{\alpha}\ln \tmfrac{6[f(x_0) - f^*]^2}{\epsilon^2}} = \tilde{\cO}\brr*{\brs{\cR\lipsqrt(\epsilon^2)/\epsilon}^2},
  \end{equation}
  where $\tilde{\cO}\brr{}$ hides logarithmic factors depending on $f(x_0) - f^*$ and $\epsilon$.
 \end{theorem}

In contrast to the projected gradient descent with scheduled stepsizes in \Cref{sec:det}, \Cref{alg:adamw} offers two substantial advantages. First, it provably works in the stochastic setting, while near-optimal stochastic guarantees for the projected gradient descent are missing. Second, \Cref{alg:adamw} does not require a complex stepsize schedule; it only requires two parameters $\alpha$ and $\cR$.

% in the case, where $\lipsqrt(s) \equiv \cM_0 + \cM_1 s$ for some possitive constants $\cM_0,\cM_1 > 0$, \Cref{ass:Lip} is equivalent to the $(\cM_0,\cM_1)$-Lipschitzness assumption considered by \citet{kovalev2026stochastic}. It is also
% Generalized Lipschitzness was previously considered by \citet{kovalev2026stochastic}

% By standard convex analysis, the subdifferential $\partial f(x)$ of the objective function $$ is nonempty at every point $x \in \Rd$.

\bibliographystyle{apalike}
\bibliography{references}

@article{kovalev2026stochastic,
  title   = {Stochastic Non-Smooth Convex Optimization with Unbounded Gradients},
  author  = {Kovalev, Dmitry},
  journal = {arXiv preprint arXiv:2605.15522},
  year    = {2026}
}

@article{alimisis2025we,
  title   = {Why do we need warm-up? a theoretical perspective},
  author  = {Alimisis, Foivos and Islamov, Rustem and Lucchi, Aurelien},
  journal = {arXiv preprint arXiv:2510.03164},
  year    = {2025}
}

@article{zhang2019gradient,
  title   = {Why gradient clipping accelerates training: A theoretical justification for adaptivity},
  author  = {Zhang, Jingzhao and He, Tianxing and Sra, Suvrit and Jadbabaie, Ali},
  journal = {arXiv preprint arXiv:1905.11881},
  year    = {2019}
}

@misc{grimmer2017convergence,
  title  = {Convergence Rates for Deterministic and Stochastic Subgradient Methods Without Lipschitz Continuity},
  author = {Grimmer, Benjamin},
  year   = {2017},
  month  = {December},
  note   = {Preprint, Optimization Online}
}

@book{nesterov2013introductory,
  title     = {Introductory lectures on convex optimization: A basic course},
  author    = {Nesterov, Yurii},
  year      = {2013},
  publisher = {Springer Science \& Business Media}
}

@book{rockafellar1997convex,
  title     = {Convex analysis},
  author    = {Rockafellar, R Tyrrell},
  volume    = {28},
  year      = {1997},
  publisher = {Princeton university press}
}

@article{zamani2025exact,
  title     = {Exact convergence rate of the last iterate in subgradient methods},
  author    = {Zamani, Moslem and Glineur, Fran{\c{c}}ois},
  journal   = {SIAM Journal on Optimization},
  volume    = {35},
  number    = {3},
  pages     = {2182--2201},
  year      = {2025},
  publisher = {SIAM}
}

@article{liu2025theoretical,
  title   = {Theoretical analysis on how learning rate warmup accelerates convergence},
  author  = {Liu, Yuxing and Ge, Yuze and Pan, Rui and Kang, An and Zhang, Tong},
  journal = {arXiv preprint arXiv:2509.07972},
  year    = {2025}
}

@article{riabinin2026does,
  title   = {Where Does Warm-Up Come From? Adaptive Scheduling for Norm-Constrained Optimizers},
  author  = {Riabinin, Artem and Veprikov, Andrey and Bolatov, Arman and Tak{\'a}{\v{c}}, Martin and Beznosikov, Aleksandr},
  journal = {arXiv preprint arXiv:2602.05813},
  year    = {2026}
}

@article{orabona2018scale,
  title     = {Scale-free online learning},
  author    = {Orabona, Francesco and P{\'a}l, D{\'a}vid},
  journal   = {Theoretical Computer Science},
  volume    = {716},
  pages     = {50--69},
  year      = {2018},
  publisher = {Elsevier}
}

@article{ahn2024adam,
  title   = {Adam with model exponential moving average is effective for nonconvex optimization},
  author  = {Ahn, Kwangjun and Cutkosky, Ashok},
  journal = {Advances in Neural Information Processing Systems},
  volume  = {37},
  pages   = {94909--94933},
  year    = {2024}
}

@article{arjevani2019oracle,
  title     = {Oracle complexity of second-order methods for smooth convex optimization},
  author    = {Arjevani, Yossi and Shamir, Ohad and Shiff, Ron},
  journal   = {Mathematical Programming},
  volume    = {178},
  number    = {1},
  pages     = {327--360},
  year      = {2019},
  publisher = {Springer}
}

\newpage

\appendix

\crefalias{section}{appendix}
\crefalias{subsection}{subappendix}
\crefalias{subsubsection}{subsubappendix}

\section{Proofs}

\proofsubsection{lem:lip}

Let $\phi(x) = \Phi(f(x) - f^*)$. Suppose that \Cref{ass:Lip} holds. Since the function $f\brr{}$ is convex on $\Rd$, it is locally Lipschitz. Hence, the functions $f\brr{}$ and $\phi\brr{}$ are differentiable almost everywhere. For arbitrary point $x \in \Rd$, where the function $f\brr{}$ is differentiable, we can obtain
\begin{align*}
  \norm{\nabla \phi(x)}
  =
  \norm{\Phi'(f(x) - f^*)\cdot \nabla f(x)}
  \at{uses \Cref{ass:Lip} and the definition in \cref{eq:Phi}}\leq
  1,
\end{align*}
where \annotate. Hence, we can prove $1$-Lipschitzness of the function $\phi\brr{}$ using Lebesgue's theorem.

Suppose that the function $\phi\brr{}$ is $1$-Lipschitz. Let $x(t) = x + t\nabla f(x)/\norm{\nabla f(x)}$, where $x \in \Rd$ is such that $\nabla f(x) \neq 0$. Then, for all $t > 0$ we can obtain the following:
\begin{align*}
  t = \norm{x(t) - x}
  \geq
  \phi(x(t)) - \phi(x)
  \geq
  \Phi(f(x) - f^* + t\norm{\nabla f(x)} ) - \Phi(f(x)-f^*).
\end{align*}
Dividing by $t$ and taking the limit $t \downarrow 0$ gives
\begin{align*}
  1 \geq \Phi'(f(x) - f^*) \cdot \norm{\nabla f(x)}
  = \norm{\nabla f(x)}/\lipsqrt([f(x) - f^*]^2),
\end{align*}
which concludes the proof.\qed

\proofsubsection{lem:f_hard}

First, it is easy to verify that the function $\Phi(s)$ is a continuous strictly increasing function. Moreover, one can verify that
\begin{align*}
  \lim_{s \to +\infty} \Phi(s)
  &=
  \medint\int_0^{+\infty}\tmfrac{\mathrm{d}t}{\lipsqrt(t^2)}
  \geq
  \medint\int_1^{+\infty}\tmfrac{\mathrm{d}t}{\lipsqrt(t^2)}
  \at{uses the concavity and positivity of the function $\lip\brr{}$}\geq
  \medint\int_1^{+\infty}\tmfrac{\mathrm{d}t}{t\lipsqrt(1)}
  =+\infty.
\end{align*}
where \annotate. Hence, its inverse $\cL\brr{}\colon \R_+ \to \R_+$ is a well-defined strictly increasing function. Furthermore, by standard calculus, we have
\begin{align*}
  \cL'(s) = 1/ \Phi'(\cL(s)) = \lipsqrt(\cL^2(s)).
\end{align*}
Hence, the function $\cL(s)$ is convex, and the function $f(x)$ in \cref{eq:f_hard} is well-defined, convex. Moreover, it is easy to observe that $f^* = 0$. 
Furthermore, let $v(x) \in \partial h_{d,\rho}(x)$ for all $x \in \Rd$. From standard results in convex analysis \citep{rockafellar1997convex} it follows that
\begin{align*}
  \partial f(x) \supset \cL'(h_{d,\rho}(x))\cdot v(x)
  =
  \lipsqrt(\cL^2(h_{d,\rho}(x)))\cdot v(x)
  \at{uses the definition in \cref{eq:f_hard}}=
  \lipsqrt([f(x) - f^*]^2)\cdot v(x),
\end{align*}
where \annotate. From the definition in \cref{eq:h_hard} it follows $\norm{v(x)} \leq 1$ for all $x \in \Rd$. Using the convexity of $f(x)$, we can conclude that $\sqn{g} \leq \lip([f(x) - f^*]^2)$ for all $x \in \Rd$ and $g \in \partial f(x)$. Hence, the function $f\brr{}$ defined in \cref{eq:f_hard} satisfies \Cref{ass:Lip}. Finally, it is easy to show that the point $x^* = -\rho(1,\dots,1)\in \Rd$ is a global solution to problem~\eqref{eq:main}.\qed

\proofsubsection{thm:lb}

Let $\rho = \cR / \sqrt{d}$ and $d = \floor{\cR/\Phi(\epsilon)}^2$. Note that $d \geq 1$ due to the fact that $\Phi(\epsilon) \leq \cR$. We choose $f\brr{}\colon \Rd \to \R$ and $\nabla f\brr{}\colon \Rd \to \Rd$ according to \cref{eq:f_hard,eq:grad_hard}. Using \Cref{lem:grad_hard} and standard reasoning \citep{nesterov2013introductory}, we can lower bound $f(x_k) - f^*$ for $k < d$ as follows: 
\begin{align*}
  f(x_k) - f^*
  &\at{uses the definition in \cref{eq:f_hard} and the fact that $h_{d,\rho}(x_k) \geq \rho$ for $k < d$}\geq
  \Phi^{-1}(\rho)
  \at{uses the definition of $\rho$ and $d$}=
  \Phi^{-1}(\cR / \floor{\cR/\Phi(\epsilon)})
  \geq \epsilon,
\end{align*}
where \annotate. Hence, to ensure $f(x_K) - f^* < \epsilon$ it must hold that $K \geq d$, which implies the following:
\begin{align*}
  K
  &\geq
  d
  \geq
  \Omega\brr{\cR^2/\Phi^2(\epsilon)}.
\end{align*}
It remains to upper-bound $\Phi(\epsilon)$. In the case $\lip(\epsilon^2) \leq 2\lip(0)$, it is easy to verify that $\Phi(\epsilon) \leq 2\epsilon / \lipsqrt(\epsilon^2)$. Otherwise, we can upper-bound $\Phi(\epsilon)$ as follows:
\begin{align*}
  \Phi(\epsilon)
  &\at{uses the definition in \cref{eq:Phi}}=
  \medint\int_0^{\epsilon}\tmfrac{\mathrm{d}t}{\lipsqrt(t^2)}
  \at{uses the concavity of the function $\lip\brr{}$}\leq
  \medint\int_0^{\epsilon}\tmfrac{\mathrm{d}t}{\brs{\lip(0) + (\lip(\epsilon^2) - \lip(0)) (t/\epsilon)^2}^{1/2}}
  \\&=
  \tmfrac{\epsilon}{\brs{\lip(\epsilon^2) - \lip(0)}^{1/2}}
  \ln\tmfrac{\lipsqrt(\epsilon^2) + \brs{\lip(\epsilon^2) - \lip(0)}^{1/2}}{\lipsqrt(0)}
  \\&\leq
  \tmfrac{2\epsilon}{\lipsqrt(\epsilon^2)}
  \ln\tmfrac{2\lipsqrt(\epsilon^2)}{\lipsqrt(0)},
\end{align*}
where \annotate.\qed

\proofsubsection{lem:GD}

Let $z_l,\dots,z_k \in \Rd$ be defined as follows:
\begin{align*}
  z_j = \argmin_{x \in Q_\cR} \brs*{H_{k,j} f(x) + \tfrac{1}{4}\sqn{x - x_j}}.
\end{align*}
Then, using the optimality of $z_j$ and the strong convexity of the objective function above, we can obtain the following inequality:
\begin{align*}
  H_{k,j} f(z_j) + \tfrac{1}{4}\sqn{z_j - x_j} + \tfrac{1}{4}\sqn{z_j - x_j}
  \leq
  H_{k,j}f(x_j).
\end{align*}
After rearranging, we get the following inequality
\begin{align*}
  \tmfrac{1}{2H_{k,j}}\sqn{z_j - x_j}
  \leq f(x_j) - f(z_j)
  \at{uses the convexity of $f\brr{}$}\leq \<\nabla f(x_j),x_j - z_j>.
\end{align*}
where \annotate.
Furthermore, we obtain the following inequality:
\begin{align*}
  \sqn{x_{j+1} - z_j}
  &\at{uses the definition in \cref{eq:GD} and the non-expansiveness of the projection}\leq
  \sqn{x_j - z_j - \eta_j\nabla f(x_j)}
  =
  \sqn{x_j - z_j}
  -2\eta_j\<\nabla f(x_j),x_j - z_j>
  +\eta_j^2\sqn{\nabla f(x_j)}
  \\&\at{uses the inequality above}\leq
  \tmfrac{H_{k,j} - \eta_j}{H_{k,j}}\sqn{z_j - x_j}
  +\eta_j^2\sqn{\nabla f(x_j)}
  \\&\at{uses the definition of $H_{k,j}$}=
  \tmfrac{H_{k,j+1}}{H_{k,j}}\sqn{z_j - x_j}
  +\eta_j^2\sqn{\nabla f(x_j)},
\end{align*}
where \annotate. After dividing by $4H_{k,j+1}$ and rearranging, we get
\begin{align*}
  f(z_j) + \tmfrac{1}{4H_{k,j}}\sqn{z_j - x_j}
  +\tmfrac{\eta_j^2}{4H_{k,j+1}}\sqn{\nabla f(x_j)}
  &\geq 
  \tmfrac{1}{4H_{k,j+1}}\sqn{x_{j+1} - z_j} + f(z_j)
  \\&\at{uses the optimality of $z_{j+1}$ above}\geq
  \tmfrac{1}{4H_{k,j+1}}\sqn{x_{j+1} - z_{j+1}} + f(z_{j+1}),
\end{align*}
where \annotate. Now, we sum these inequalities for $j = l,\dots,k-1$ and obtain
\begin{align*}
  \tmfrac{1}{4\eta_k}\sqn{x_k - z_k} + f(z_k) &\leq f(z_l) + \tmfrac{1}{4H_{k,l}}\sqn{z_l - x_l}
  +\msum_{j=l}^{k-1}\tmfrac{\eta_j^2}{4H_{k,j+1}}\sqn{\nabla f(x_j)}
  \\&\at{uses the optimality of $z_l$ above}\leq
  f(x) + \tmfrac{1}{4H_{k,l}}\sqn{x_l - x}
  +\msum_{j=l}^{k-1}\tmfrac{\eta_j^2}{4H_{k,j+1}}\sqn{\nabla f(x_j)},
\end{align*}
where \annotate. After rearranging, we get
\begin{align*}
  \mi{10}\tmfrac{1}{4H_{k,l}}\sqn{x_l - x}
  +\msum_{j=l}^{k-1}\tmfrac{\eta_j^2}{4H_{k,j+1}}\sqn{\nabla f(x_j)}
  \geq
  \tmfrac{1}{4\eta_k}\sqn{x_k - z_k} + f(z_k) - f(x)
  \\&\at{uses the convexity of $f\brr{}$}\geq
  \tmfrac{1}{4\eta_k}\sqn{x_k - z_k} + f(x_k) + \<\nabla f(x_k),z_k - x_k> - f(x)
  \\&\at{uses Young's inequality}\geq
  f(x_k) - f(x) - \eta_k\sqn{\nabla f(x_k)},
\end{align*}
where \annotate.\qed

\proofsubsection{thm:GD}

We prove that $f(x_k) - f^* \leq F_k$, $0 < \alpha_k < \frac{1}{32}$ and $0 < F_{k+1} < F_k$ by induction over $k \in \N_0$. The base case $k=0$ can be proved directly. Now, suppose that $f(x_j) - f^* \leq F_j$, $0 < \alpha_j < \frac{1}{32}$, and $0 < F_{j+1} < F_j$ for all $0 \leq j < k$. Using the concavity and monotonicity of the function $\lip\brr{}$, we get
\begin{align*}
  \tmfrac{F_k}{\lipsqrt(F_k^2)} \leq \tmfrac{F_0}{\lipsqrt(F_0^2)} = \tmfrac{3\Delta_0}{\lipsqrt(9\Delta_0^2)} \leq \tmfrac{3\Delta_0}{\lipsqrt(\Delta_0^2)} \leq 3\cR.
\end{align*}
Hence, we can bound $\alpha_k$ as follows:
\begin{align*}
  0 < \alpha_k
  \at{use the definitions in \cref{eq:stepsize}}
  \leq \tmfrac{F_k^2}{8\cR^2L_k \lip(F_k^2)}
  \at{uses the inequality above}\leq
  \tmfrac{9}{8L_k}
  \at{use the definitions in \cref{eq:stepsize}}< \tfrac{1}{32}
\end{align*}
where \annotate. Hence, using the definition in \cref{eq:stepsize} we conclude that
\begin{align*}
  0 < F_{k+1} = (1-\alpha_k)F_k < F_k.
\end{align*}
Next, we can upper-bound $f(x_k) - f^*$ as follows:
\begin{align*}
  f(x_k) - f^*
  &\at{uses the convexity of $f\brr{}$ and the induction hypothesis}\leq
  F_{k-1} + \<\nabla f(x_k),x_k - x_{k-1}>
  \\&\at{uses \Cref{ass:Lip} and \cref{eq:GD}}\leq
  F_{k-1} + \eta_{k-1}\lipsqrt(F_{k-1}^2)\lipsqrt([f(x_k) - f^*]^2)
  \\&\at{uses the definition in \cref{eq:stepsize}}\leq
  F_{k-1} + \tmfrac{F_{k-1}\lipsqrt([f(x_k) - f^*]^2)}{L_{k-1}\lipsqrt(F_{k-1}^2)}
\end{align*}
where \annotate. In the case $f(x_k) - f^* > F_{k-1}$, we can obtain
\begin{align*}
  f(x_k) - f^* \leq
  F_{k-1} + \tmfrac{F_{k-1}\lipsqrt([f(x_k) - f^*]^2)}{L_{k-1}\lipsqrt(F_{k-1}^2)}
  \at{uses the concavity of the function $\lip\brr{}$}\leq
   F_{k-1} + \tmfrac{f(x_k) - f^*}{L_{k-1}},
\end{align*}
where \annotate, and the same inequality obviously holds in the case $f(x_k) - f^* \leq F_{k-1}$. Hence, using the definition in \cref{eq:stepsize}, we can obtain
\begin{align*}
  f(x_k) - f^*
  \leq
  \tfrac{64}{63}F_{k-1}
  \at{uses the definition in \cref{eq:stepsize}}=
  \tmfrac{64}{63(1-\alpha_{k-1})}F_k
  \at{uses the induction hypothesis}\leq
  \tmfrac{64\cdot 32}{63\cdot 31}F_k < 2F_k,
\end{align*}
where \annotate. Hence, using \Cref{ass:Lip}, we can upper-bound $\norm{\nabla f(x_k)}$ as follows:
\begin{align*}
  \norm{\nabla f(x_k)}
  \leq
  \lipsqrt([f(x_k) - f^*]^2)
  \leq
  \lipsqrt(4F_k^2)
  \at{uses the concavity and positivity of the function $\lip\brr{}$}\leq
  2\lipsqrt(F_k^2),
\end{align*}
where \annotate. Now, we define $l \in \N_0$ as follows:
\begin{align*}
  l = \min\brf{j \in \brf{0,\dots,k} : F_j \leq \tfrac{3}{2}F_k}.
\end{align*}
For all $l \leq j \leq k$, we get $F_k \leq F_j \leq \frac{3}{2}F_k$. Moreover, for all $l \leq j \leq k$, we can lower-bound $\eta_j$ as
\begin{align*}
  \eta_j
  \at{use the definitions in \cref{eq:stepsize}}= 
  \tmfrac{F_j^2}{L_jF_j\lip(F_j^2)}
  \at{uses the concavity and positivity of the function $\lip\brr{}$}\geq
  \tmfrac{F_k^2}{L_jF_j\lip(F_k^2)}
  \geq
  \tmfrac{2F_k}{3L_j\lip(F_k^2)}
  \at{use the definitions in \cref{eq:stepsize}}\geq 
  \tmfrac{2F_k}{3L_k\lip(F_k^2)},
\end{align*}
and lower-bound $H_{k,j}$ as
\begin{align*}
  H_{k,j}
  =
  \msum_{i = j}^k \eta_i
  \geq
  \tmfrac{2F_k(k-j+1)}{3L_k\lip(F_k^2)},
\end{align*}
where \annotate. Using \Cref{lem:GD}, we get
\begin{align*}
  f(x_k) - f(x)
  &\leq
  \tmfrac{\sqn{x_l - x}}{4H_{k,l}} + \msum_{j=l}^{k-1} \tmfrac{\eta_j^2\sqn{\nabla f(x_j)}}{4H_{k,j+1}} + \eta_k\sqn{\nabla f(x_k)}
  \\&\at{uses \Cref{ass:Lip}, the induction hypothesis, and the previously obtained upper bound on $\norm{\nabla f(x_k)}$}\leq
  \tmfrac{\sqn{x_l - x}}{4H_{k,l}} + \msum_{j=l}^{k-1} \tmfrac{\eta_j^2\lip(F_j^2)}{4H_{k,j+1}} + 4\eta_k\lip(F_k^2)
  \\&\at{uses the inequalities above}\leq
  \tmfrac{\sqn{x_l - x}}{4H_{k,l}} + \msum_{j=l}^{k-1} \tmfrac{3L_k\eta_j^2\lip(F_j^2)\lip(F_k^2)}{8(k-j)F_k} + 4\eta_k\lip(F_k^2)
  \\&\at{use the definitions in \cref{eq:stepsize}}=
  \tmfrac{\sqn{x_l - x}}{4H_{k,l}} + \msum_{j=l}^{k-1} \tmfrac{3L_k\lip(F_k^2)}{8(k-j)L_j^2F_k}\cdot\tmfrac{F_j^2}{\lip(F_j^2)} + \tmfrac{4F_k}{L_k}
  \\&\at{uses the monotonicity of the function $\lip\brr{}$}\leq
  \tmfrac{\sqn{x_l - x}}{4H_{k,l}} + \msum_{j=l}^{k-1} \tmfrac{3L_k\lip(F_k^2)}{8(k-j)L_j^2F_k}\cdot\tmfrac{9F_k^2}{4\lip(F_k^2)} + \tmfrac{4F_k}{L_k}
  \\&=
  \tmfrac{\sqn{x_l - x}}{4H_{k,l}} + \brs*{\msum_{j=l}^{k-1} \tmfrac{27L_k}{32(k-j)L_j^2} + \tmfrac{4}{L_k}}F_k
  \\&\at{use the definitions in \cref{eq:stepsize}}\leq
  \tmfrac{\sqn{x_l - x}}{4H_{k,l}} + \tfrac{43}{256}F_k,
\end{align*}
where \annotate. Furthermore, in the case $l = 0$, we can choose $x = x_0$ and upper-bound $f(x_k) - f^*$ as follows:
\begin{align*}
  f(x_k) - f^* \leq f(x_0) - f^* + \tmfrac{43}{256}F_k \leq \tfrac{1}{3}F_0 + \tmfrac{43}{256}F_k
  \leq 
  \tfrac{1}{2}F_k + \tfrac{43}{256}F_k \leq F_k.
\end{align*}
Finally, in the case $l \neq 0$, we have $F_{l-1} > \frac{3}{2}F_k$, which implies
\begin{align*}
  F_k < \tfrac{2}{3}F_{l-1} = \tmfrac{2}{3(1-\alpha_{l-1})}F_l \leq \tfrac{64}{93}F_l.
\end{align*}
Therefore, we get $l < k$. Furthermore, we can lower-bound $H_{k,l}$ as follows:
\begin{align*}
  H_{k,l}
  \geq
  \msum_{j=l}^{k-1} \eta_j
  &\at{uses the definitions in \cref{eq:stepsize}}=
  \msum_{j=l}^{k-1} \tmfrac{8\cR^2 (1-\alpha_j)(F_j - F_{j+1})}{F_jF_{j+1}}
  \\&\at{use the previously obtained inequalities}\geq
  \msum_{j=l}^{k-1} \tfrac{31}{4}\cR^2\brs*{\tmfrac{1}{F_{j+1}} - \tmfrac{1}{F_j}}
  =
  \tfrac{31}{4}\cR^2\brs*{\tmfrac{1}{F_k} - \tmfrac{1}{F_l}}
  \at{use the previously obtained inequalities}\geq
  \tmfrac{29\cR^2}{12F_k},
\end{align*}
where \annotate. Hence, we obtain
\begin{align*}
  f(x_k) - f^* &\leq \tmfrac{\sqn{x_l - x^*}}{4H_{k,l}} + \tmfrac{43}{256}F_k
  \\&\at{uses the fact that $x_l,x^* \in Q_\cR$}\leq
  \tmfrac{\cR^2}{H_{k,l}} + \tmfrac{43}{256}F_k
  \\&\at{uses the inequality above}\leq
  \brs*{\tfrac{12}{29} + \tfrac{43}{256}}F_k \leq F_k,
\end{align*}
where \annotate.\qed

\proofsubsection{cor:GD}

Since $F_k \geq F_{k+1}$, there exists the limit $F = \lim_{k\to \infty} F_k \geq 0$. Suppose that $F>0$. Then we obtain the inequality $\alpha_k \geq \frac{F^2}{8\cR^2L_k\lip(F_0^2)}$, which implies
\begin{align*}
  -\infty < \ln F = \ln F_0 +  \msum_{j = 0}^{\infty}\ln(1-\alpha_j) \leq \ln F_0 -\msum_{j = 0}^{\infty} \tmfrac{F^2}{8\cR^2L_j\lip(F_0^2)} = -\infty,
\end{align*}
which is not possible. Hence, $F = 0$, and we can define $K = \min\brf{k : F_k \leq \epsilon}$. In the case $K = 0$, the desired result is obvious. Otherwise, we can obtain the following:
\begin{align*}
  % \brr*{\tmfrac{\cR\lipsqrt(\epsilon^2)}{\epsilon}}^2
  \cR^2\medint\int_{\epsilon^2}^{F_0^2} \tmfrac{\lip(s)\mathrm{d}s}{s^2}
  \at{uses the concavity and positivity of the function $\lip\brr{}$}\leq 
  \cR^2\medint\int_{\epsilon^2}^{F_0^2} \tmfrac{\lip(\epsilon^2)\mathrm{d}s}{\epsilon^2 s}
  \leq
  \tilde{\cO}\brr*{\brs{\cR\lipsqrt(\epsilon^2)/\epsilon}^2}.
\end{align*}
where \annotate.
Consequently, we can obtain the following inequality:
\begin{align*}
  \tilde{\cO}\brr*{\brs{\cR\lipsqrt(\epsilon^2)/\epsilon}^2}
  &\geq
  \cR^2\medint\int_{\epsilon^2}^{F_0^2} \tmfrac{\lip(s)\mathrm{d}s}{s^2}
  \geq
  \cR^2\msum_{k=1}^{K-1}\medint\int_{F_k^2}^{F_{k-1}^2} \tmfrac{\lip(s)\mathrm{d}s}{s^2}
  \\&\at{uses the concavity and positivity of the function $\lip\brr{}$}\geq 
  \cR^2\msum_{k=1}^{K-1}\medint\int_{F_k^2}^{F_{k-1}^2} \tmfrac{\lip(F_{k-1}^2)\mathrm{d}s}{sF_{k-1}^2}
  =
  2\cR^2\msum_{k=1}^{K-1}\tmfrac{\lip(F_{k-1}^2)}{F_{k-1}^2}\ln\tmfrac{F_{k-1}}{F_k}
  \\&\at{use the definitions in \cref{eq:stepsize}}=
  -2\cR^2\msum_{k=1}^{K-1}\tmfrac{\lip(F_{k-1}^2)}{F_{k-1}^2}\ln(1-\alpha_{k-1})
  \geq
  2\cR^2\msum_{k=1}^{K-1}\tmfrac{\alpha_{k-1}\lip(F_{k-1}^2)}{F_{k-1}^2}
  \\&\at{use the definitions in \cref{eq:stepsize}}=
  \msum_{k=1}^{K-1}\tmfrac{1}{4L_{k-1}} \geq \tmfrac{K-1}{64(4 + \ln K)},
\end{align*}
where \annotate, which implies $K \leq \tilde{\cO}\brr*{\brs{\cR\lipsqrt(\epsilon^2)/\epsilon}^2}$.\qed

\proofsubsection{lem:adamw}

 Let $z_{k+1} = -\cR\clip{m_k / \sqrt{s_k}}$. Then, from \Cref{line:x} it follows that $x_{k+1} = (1-\alpha)x_k + \alpha z_{k+1}$. Furthermore, we can upper-bound $(1-\alpha)[f(x_{k+1}) - f^*]$ as follows:
 \begin{align*}
  \mi{2}(1-\alpha)[f(x_{k+1}) - f^*]
  \at{use the convexity of the function $f\brr{}$}\leq
  (1-\alpha)\brs{f(x_k) - f^* + \<\nabla f(x_{k+1}),x_{k+1} - x_k>}
  \\&=
  (1-\alpha)[f(x_k) - f^*] + \alpha\<\nabla f(x_{k+1}),z_{k+1} - x^* + x^* - x_{k+1}>
  \\&\at{use the convexity of the function $f\brr{}$}\leq
  (1-\alpha)[f(x_k) - f^*] - \alpha[f(x_{k+1}) - f^*] + \alpha\<\nabla f(x_{k+1}),z_{k+1} - x^*>,
 \end{align*}
 where \annotate. After rearranging, we obtain
 \begin{align*}
  f(x_{k+1}) - f^* \leq (1-\alpha)[f(x_k) - f^*] + \alpha\<\nabla f(x_{k+1}),z_{k+1} - x^*>.
 \end{align*}
 After combining these inequalities for $k = 0,\dots,K-1$, we get
\begin{align*}
  f(x_K) - f^*
  &\leq
  (1-\alpha)^K[f(x_0) - f^*] + \msum_{k=1}^K\alpha(1-\alpha)^{K-k}\<\nabla f(x_k),z_k - x^*>
  \\&=
  (1-\alpha)^K[f(x_0) - f^*] + \msum_{k=1}^K\alpha(1-\alpha)^{K-k}\<\nabla f(x_k) - g_k,z_k - x^*>
  \\&+
  \msum_{k=1}^K\alpha(1-\alpha)^{K-k}\<g_k,z_k - x^*>.
\end{align*}
One can verify that the update rule for $z_k$ is equivalent to an iteration of the SOLO-FTRL algorithm \citep{orabona2018scale}. Hence, we can use Lemma~8 of \citet{ahn2024adam} and obtain
\begin{align*}
  f(x_K) - f^*
  &\leq
  (1-\alpha)^K[f(x_0) - f^*] + \msum_{k=1}^K\alpha(1-\alpha)^{K-k}\<\nabla f(x_k) - g_k,z_k - x^*>
  \\&+
  4\alpha\cR\sqrt{\msum_{k=1}^K(1-\alpha)^{2(K-k)}\sqn{g_k}}.
\end{align*}
Next, we square these inequality and take expectation, which gives the following:
\begin{align*}
  \mi{5}\E{(f(x_K) - f^*)^2}
  \leq
  3(1-\alpha)^{2K}[f(x_0) - f^*]^2 +  48\alpha^2\cR^2\msum_{k=1}^K(1-\alpha)^{2(K-k)}\E{\sqn{g_k}}
  \\&
  +3\E*{\brr*{\msum_{k=1}^K\alpha(1-\alpha)^{K-k}\<\nabla f(x_k) - g_k,z_k - x^*>}^2}
  \\&\at{uses the unbiasedness of $g_k$ and the independence of samples $\xi_k$ on \Cref{line:g}}=
  3(1-\alpha)^{2K}[f(x_0) - f^*]^2 +  48\alpha^2\cR^2\msum_{k=1}^K(1-\alpha)^{2(K-k)}\E{\sqn{g_k}}
  \\&
  +3\E*{\msum_{k=1}^K\alpha^2(1-\alpha)^{2(K-k)}\<\nabla f(x_k) - g_k,z_k - x^*>^2}
  \\&\at{uses the fact that $\norm{z_k - x^*} \leq 2\cR$}\leq
  3(1-\alpha)^{2K}[f(x_0) - f^*]^2 +  48\alpha^2\cR^2\msum_{k=1}^K(1-\alpha)^{2(K-k)}\E{\sqn{g_k}}
  \\&
  +12\alpha^2\cR^2\msum_{k=1}^K(1-\alpha)^{2(K-k)}\E*{\sqn{\nabla f(x_k) - g_k}}
  \\&\at{uses the unbiasedness of $g_k$}\leq
  3(1-\alpha)^{2K}[f(x_0) - f^*]^2 +  60\alpha^2\cR^2\msum_{k=1}^K(1-\alpha)^{2(K-k)}\E{\sqn{g_k}}
  \\&\at{uses \Cref{ass:Lip_stoch} and the concavity of the function $\lip\brr{}$}\leq
   3(1-\alpha)^{2K}[f(x_0) - f^*]^2 +  60\alpha^2\cR^2\msum_{k=1}^K(1-\alpha)^{2(K-k)}\lip(\E{(f(x_k) - f^*)^2}),
\end{align*}
where \annotate.\qed

\proofsubsection{lem:M01}

in the case $s \leq \epsilon$, the desired inequality is implied by the monotonicity of the function $\lip\brr{}$. Otherwise, using the positivity and concavity of the function $\lip\brr{}$, we can obtain the inequality $\lip(s^2)/s^2 \leq \lip(\epsilon^2)/\epsilon^2$, which implies the desired upper bound.\qed

\proofsubsection{thm:adamw}

 From \Cref{lem:adamw,lem:M01} it follows that:
\begin{align*}
  \Delta_K^2 \leq 3(1-\alpha)^{2K}\Delta_0^2 + [8\alpha\cR]^2\msum_{k=1}^K (1-\alpha)^{2(K-k)}(\cM_0^2 + \cM_1^2\Delta_k^2).
\end{align*}
Let the sequence $\hDelta_k$ be defined recursively as follows:
\begin{align*}
  \hDelta_0^2 = 3\Delta_0^2,\quad
  \hDelta_{k+1}^2 =\tmfrac{(1-\alpha)^2\hDelta_k^2 + [8\alpha\cR\cM_0]^2}{(1 - [8\alpha\cR\cM_1]^2)}.
\end{align*}
Using this definition, we can show that 
\begin{align*}
  (1-\alpha)^{2(K-k-1)}\hDelta_{k+1}^2
  &=
  (1-\alpha)^{2(K-k)}\hDelta_k^2
  \\&
  +[8\alpha\cR\cM_1]^2(1-\alpha)^{2(K-k-1)}\hDelta_{k+1}^2 + (1-\alpha)^{2(K-k-1)}[8\alpha\cR\cM_0]^2.
\end{align*}
Now, we sum these inequalities for $k = 0,\dots,K-1$ and obtain the following:
\begin{align*}
  \hDelta_K^2 = (1-\alpha)^{2K}\hDelta_0^2 + [8\alpha\cR]^2\msum_{k=1}^K
  (1-\alpha)^{2(K-k)}[\cM_0^2 + \cM_1^2\hDelta_k^2].
\end{align*}
Next, we will prove that $\Delta_K^2 \leq \hDelta_K^2$ by induction over $K$. The base case $K=0$ is obvious. For $K \in \N$, we can upper-bound $\Delta_K^2$ as follows:
\begin{align*}
  \mi{3}[1 - [8\alpha\cR\cM_1]^2]\Delta_K^2
  \\&\leq
  3(1-\alpha)^{2K}\Delta_0^2 + [8\alpha\cR]^2\msum_{k=1}^K (1-\alpha)^{2(K-k)}\cM_0^2
  +[8\alpha\cR\cM_1]^2\msum_{k=1}^{K-1} (1-\alpha)^{2(K-k)}\Delta_k^2
  \\&\at{uses the induction hypothesis}\leq
  (1-\alpha)^{2K}\hDelta_0^2 + [8\alpha\cR]^2\msum_{k=1}^K (1-\alpha)^{2(K-k)}\cM_0^2
  +[8\alpha\cR\cM_1]^2\msum_{k=1}^{K-1} (1-\alpha)^{2(K-k)}\hDelta_k^2
  \\&\at{uses the previously obtained relation for $\hDelta_K$}\leq
  [1 - [8\alpha\cR\cM_1]^2]\hDelta_K^2,
\end{align*}
where \annotate. Using the fact that $8\alpha\cR\cM_1 < 1$, we obtain the desired inequality $\Delta_K^2 \leq \hDelta_K^2$. Now, we can upper bound $\hDelta_{k+1}^2$ as follows:
\begin{align*}
 \hDelta_{k+1}^2
 &=
 \tmfrac{(1-\alpha)^2\hDelta_k^2 + [8\alpha\cR\cM_0]^2}{(1 - [8\alpha\cR\cM_1]^2)}
 \at{uses the definition of $\cM_1$ and $\alpha$}\leq
 (1-\alpha)\hDelta_k^2 + \tfrac{2}{7}\alpha\epsilon^2,
\end{align*}
where \annotate. This implies
\begin{align*}
  \Delta_K^2 \leq \hDelta_K^2 \leq 3(1-\alpha)^K\Delta_0^2 + \tfrac{2}{7}\epsilon^2
  \leq \tfrac{1}{2}\epsilon^2 + \tfrac{2}{7}\epsilon^2 \leq \epsilon^2,
\end{align*}
which concludes the proof.\qed

\end{document}